\documentclass[12pt,draftcls,onecolumn]{IEEEtran}

\usepackage{amsmath}
\usepackage{amssymb}
\usepackage{wasysym}

\newcommand{\cvd}{ {\APLbox}}

\newcommand{\Range}{\mathop{\rm Range}}

\newcommand{\tr}{\mathop{\rm tr}}
\newcommand{\Rgamma}{\mathop{\rm Range} \Gamma}

\newcommand{\C}{\mathbb{C}}

\newcommand{\ejth}{{\rm e}^{{\rm j}\vartheta}}

\newcommand{\Hermitian}{\mathfrak{H}}

\newcommand{\DD}{D}

\newcommand{\Lambdaopt}{\Lambda_\circ}

\newcommand{\eajt}{{\rm e}^{j\vartheta}}

\usepackage{graphicx}
\usepackage{amssymb}
\usepackage{epstopdf}
\usepackage{amsfonts}
\usepackage{mathrsfs}
\usepackage{marvosym}
\usepackage[pdftex]{color}
\usepackage{pdfcolmk}

\usepackage[dvips]{epsfig} \DeclareGraphicsExtensions{.ps,.eps}
\renewcommand{\baselinestretch}{1.6}

\newfont{\BB}{msbm10}
\newfont{\bb}{msbm8}

\def\KLD{\mbox{\BB S}}

\def\C{\mbox{\BB C}}

\def\T{\mbox{\BB T}}
\def\J{\mbox{\BB J}}

\def\tr{{\rm tr}\,}

\newcommand{\beq}{\begin{equation}}
\newcommand{\eeq}{\end{equation}}
\newcommand{\bea}{\begin{eqnarray}}
\newcommand{\eea}{\end{eqnarray}}

\newcommand{\bsea}{\begin{subeqnarray}}
\newcommand{\esea}{\end{subeqnarray}}
\newcommand{\nn}{\nonumber}

\newcommand{\rst}{\setcounter{acount}{0}}
\newcommand{\stoa}{\setcounter{acount}{1}}
\newcommand{\stob}{\addtocounter{equation}{-1} \setcounter{acount}{2}}

\def\bmat{\left[ \begin{array}}
\def\emat{\end{array} \right]}

\newcounter{acount}
\renewcommand{\theequation}{\arabic{equation}\alph{acount}}

\newtheorem{theo}{Theorem}[section]

\newtheorem{lem}{Lemma}[section]
\newtheorem{cor}{Corollary}[section]

\newtheorem{prop}{Proposition}[section]
\newtheorem{rem}{Remark}[section]
\newtheorem{problem}{Problem}[section]

\definecolor{red}{cmyk}{0.1,0.8,0,0.1}

\begin{document}
\title{  On the convergence of an efficient algorithm for Kullback-Leibler approximation  of spectral densities}
\author{Augusto~Ferrante, Federico~Ramponi, and~Francesco~Ticozzi\thanks{Work partially supported by the Italian Ministry for Education and Resarch (MIUR) under PRIN grant
``Identification and Robust Control of Industrial Systems".}
\thanks{A. Ferrante and F. Ticozzi are  with the
Dipartimento di Ingegneria dell'Informazione, Universit\`a di Padova,
via Gradenigo 6/B, 35131
Padova, Italy. {\tt\small  augusto@dei.unipd.it},  {\tt\small  ticozzi@dei.unipd.it}}
\thanks{F. Ramponi is with   the
Institut f\"ur Automatik, ETH Z\"urich,
Physikstrasse 3, 8092 Z\"urich, Switzerland. {\tt\small
ramponif@control.ee.ethz.ch}}}

\markboth{DRAFT}{Shell \MakeLowercase{\textit{et al.}}: Bare Demo of IEEEtran.cls for Journals}

\maketitle

\begin{abstract}
This paper   deals with  a method for the approximation of a
spectral density function among the solutions of a generalized moment problem
\`a la Byrnes/Georgiou/Lindquist.
The approximation is pursued with respect to the Kullback-Leibler
pseudo-distance, which gives rise to a convex optimization problem.
After developing the variational analysis,
we discuss the properties of an   efficient 
algorithm for the solution of the   corresponding  dual problem, 
based on the iteration of   a nonlinear map  in a bounded subset of the dual space.
Our main result is the proof of  local convergence of the
latter, established as a consequence of the Central Manifold Theorem.
 Supported by numerical evidence,  we conjecture that, in the mentioned bounded set, the convergence is actually {\em global}.
\end{abstract}

\section{Introduction}

During the last decade a broad research program on the interplay between {\em (generalized) moment problems} and {\em analytic interpolation problems} with complexity constraints, {\em robust control}, {\em approximation} and {\em estimation} of spectral density functions has been carried over by C. I. Byrnes, T. Georgiou, and A. Lindquist and their co-authors and epigones
\cite{BYRNES_GUSEV_LINDQUIST_CONVEX,GEORGIOU_INTERPOLATION,BYRNES_GEORGIOU_LINDQUIST_THREE,BYRNES_GUSEV_LINDQUIST_FROMFINITE,BYRNES_GEORGIOU_LINDQUIST_MEGRETSKI,BYRNES_GEORGIOU_LINDQUIST_GENERALIZEDCRITERION,GEORGIOU_SELECTIVEHARMONIC,BYRNES_ENQVIST_LINDQUIST_IDENTIFIABILITY,BYRNES_ENQVIST_LINDQUIST_CEPSTRAL,NAGAMUNE_ROBUSTSOLVER,BLOMQVIST_LINDQUIST_NAGAMUNE_MATRIXVALUED,BYRNES_LINDQUIST_THEGENERALIZED,GEORGIOU_MAXIMUMENTROPY,GEORGIOU_RELATIVEENTROPY,GEORGIOU_RIEMANNIAN,GEORGIOU_LINDQUIST_KULLBACKLEIBLER,PAVON_FERRANTE_ONGEORGIOULINDQUIST,BYRNES_LINDQUIST_IMPORTANTMOMENTS,FERRANTE_PAVON_RAMPONI_CONSTRAINEDHELLINGER,FERRANTE_PAVON_RAMPONI_HELLINGERVSKULLBACK,RAMPONI_FERRANTE_PAVON_GLOBALLYCONVERGENT,ENQVIST_HOMOTOPYAPPROACH,GEORGIOU_DISTANCESBETWEEN,LINDQUIST_PREDICTIONERROR}.
Moment problems have a long history and have been   at the heart of  many mathematical and engineering problems in the past century, see, e.g., \cite{Shohat_Tamarkin,Akhiezer}
and the references therein.   Only with recent developments of the above-mentioned research program, however, the parametrization of solutions in the presence of additional constraints on the complexity have been satisfactorily addressed \cite{BYRNES_LINDQUIST_THEGENERALIZED}. This result, that has been possible thanks to a    suitable variational formulation, is of  key interest  in control engineering. In fact,  the well-known relation between moment problems and Nevanlinna-Pick interpolation problems   allows for solutions of  ${\cal H}_\infty$ control problems    that include  a bound on the complexity of the controller, which is of paramount practical importance \cite{BLOMQVIST_LINDQUIST_NAGAMUNE_MATRIXVALUED,BYRNES_GEORGIOU_LINDQUIST_MEGRETSKI}. Similar considerations   apply to  the {\em covariance extension} problem \cite{BYRNES_GUSEV_LINDQUIST_CONVEX,GEORGIOU_LINDQUIST_KULLBACKLEIBLER}.
Among the other applications, we also mention signal and image processing
\cite{GEORGIOU_MAXIMUMENTROPY,GEORGIOU_RELATIVEENTROPY,GEORGIOU_DISTANCESBETWEEN,GEORGIOU_RIEMANNIAN},
and Biomedical Engineering \cite{NASIRI_EBBINI_GEORGIOU_NONINVASIVEESTIMATION}.
These applications are based on a {\em spectral estimation}
 procedure that hinges on {\em optimal approximation} of spectral densities with linear {\em integral constraints} that may be viewed as constraints on a finite number of moments of the spectrum.  The linear integral constraints represent a knowledge on the steady-state state covariance of a bank of filters that is designed in order to estimate the unknown spectral density $\Phi$, while the to-be-approximated spectral density represents a prior knowledge on  $\Phi$, see Section \ref{sec:BGM} for more details.
As discussed in 
\cite{BYRNES_GEORGIOU_LINDQUIST_THREE,GEORGIOU_LINDQUIST_KULLBACKLEIBLER,PAVON_FERRANTE_ONGEORGIOULINDQUIST}, this optimal approximation leads to a tunable  spectral estimation algorithm that provides high resolution estimates in prescribed frequency bands even in presence of a short record of observed data. 
An important feature of the above mentioned optimal approximation method is that the primal optimization problem can be solved in closed form and, as long as the prior spectral density is rational, yields a rational solution 
with an {\em a priori} bound on the complexity (McMillan degree).

The numerical challenge, in practical applications, lays with the {\em dual problem}. In fact, the dual variable is an Hermitian matrix and, as discussed in \cite{GEORGIOU_LINDQUIST_KULLBACKLEIBLER}, the reparametrization in vector form  may lead to a loss of convexity.
Moreover, the dual functional and its
  gradient 
tend to infinity at the boundary, leading to serious numerical  difficulties in practical implementations.
Indeed, any {\em gradient-based} numerical method is severely affected by heavy computational burden, due to a large number of back-stepping iterations. These arise as consequence of the above-mentioned behavior of the gradient in the vicinity of the boundary. 
In order to avoid this computational slowdown, a nonlinear matricial iteration 
has  
been proposed in \cite{PAVON_FERRANTE_ONGEORGIOULINDQUIST}: It exhibits surprisingly good performance
and it does not involve either back-stepping or the computation and inversion of the Hessian. 
Proof of convergence, however, has so far been a challenging open problem.
This problem is eventually successfully addressed in this paper: Our main contribution is the proof that this iteration is locally asymptotically convergent to the manifold of solutions of the dual problem. Moreover, we analyze many other aspects of this matricial iteration and its dynamical properties in connection with the dual problem.    We finally show that, by resorting to a spectral factorization method, the proposed iteration may be implemented in an effective way.
  In fact, each iteration of the algorithm 
only requires the solutions of a Riccati equation and of a Lyapunov equation, for which robust and efficient algorithms are available.

The paper is organized as follows. In Section \ref{sec:BGM} we give a proper mathematical statement of the problem,   and proceed by recalling some relevant facts from the literature as well as establishing some preliminary results.
Section \ref{sec:POC} contains our main result: Local convergence of the matricial iteration. The proof is rather   articulated and involves many different tools from linear and non-linear systems theory, including the {\em Central Manifold Theorem}.  For this reason, this section is divided into many subsections that provide a roadmap of the   various parts of the proof.
As a byproduct, we also obtain many relevant results on the iteration and its linearization. 
In Section \ref{sec:NI} we describe how  the proposed iteration may be implemented in an effective numerical way.
Section \ref{sec:EFS}   illustrates  some results obtained from simulations and a conjecture.  Final remarks, conclusions and future perspective are presented in Section \ref{sec:C}.

{\em Notation.} 
We denote by $\Hermitian_n$ the set of Hermitian matrices of dimension $n$.
Given a complex matrix $A$, $A^*$ denotes the transpose conjugate of $A$,
while, for a matrix valued function $\chi(z)$ in the complex variable $z$,
$\chi^*(z)$ denotes the analytic continuation of the function that for $|z|=1$ equals the transpose conjugate of $\chi(z)$. Thus, for a matrix-valued rational function $\chi(z)=H(zI-F)^{-1}G+J$, we have $\chi^*(z)=G^*(z^{-1}I-F^*)^{-1}H^*+J^*$.
We denote by $\T$ the unit circle in the complex plane $\C$ and by $C(\T)$  the set of  complex-valued continuous functions on $\T$.
$C_+(\T)$ denotes the subset of $C(\T)$ whose elements are real-valued positive  functions. Elements in   $C_+(\T)$  will be thought of as {\em spectral densities}. 

\section{Problem formulation and background material}\label{sec:BGM}

Consider the rational transfer function
$$G(z)=(zI-A)^{-1}B,\qquad  A\in \C^{n\times n}, B\in\C^{n\times 1},$$
of the system
$$
x(t+1)=Ax(t)+By(t),$$
 where
$A$ is a stability matrix, i.e. has all its eigenvalues in the open  unit
disc, and $(A,B)$ is a {\em reachable pair}. 

The transfer function $G$ models a bank of filters fed by a stationary process $y(t)$ of unknown spectral density $\Phi(z)$.
We assume that we know (or that we can reliably estimate) the steady-state covariance $\Sigma$ of the state $x$ of the filter.
Based on $\Sigma$ and on an {\em a priori} information in the form of a prior spectral density $\Psi(z)$, we want to estimate the spectral density $\Phi(z).$

We will consider the Kullback-Leibler index as a measure of the difference
between spectral densities $\Psi$ and $\Phi$ in $C_+(\T)$:
$$\KLD(\Psi\|\Phi)=\int
\Psi\log\left(\frac{\Psi}{\Phi}\right)=\int_{-\pi}^{\pi}\Psi(\eajt)
\log\left(\frac{\Psi(\eajt)}{\Phi(\eajt)}\right)\,\frac
{d\vartheta}{2\pi}.
$$
The above notation,
where integration takes place on the unit circle and with respect to the normalized Lebesgue measure,
is used throughout the whole paper.

As in \cite{GEORGIOU_LINDQUIST_KULLBACKLEIBLER}, we consider
the following
\begin{problem}\label{AP} (Approximation  problem) Let $\Psi\in
C_+(\T)$, and let $\Sigma\in
\C^{n\times n}$  satisfy $\Sigma=\Sigma^*>0$. Find
$\hat{\Phi}$ that  solves
\begin{eqnarray} \label{1}&&{\rm minimize} \quad
\KLD(\Psi\|\Phi)\\\label{2}&&{\rm over}\quad\left\{\Phi\in C_+(\T)\; |
\int G\Phi  G^*=\Sigma\right\}.
\end{eqnarray}
\end{problem}

\begin{rem}
\label{REMARK_KL_POSITIVITY}
Notice that in order to guarantee that the Kullback-Leibler pseudo-distance is greater than or equal to zero we need that  the zeroth-order moment of its arguments
is the same, i.e. if $\int\Psi = \int\Phi$ then $\KLD(\Psi||\Phi) \geq 0$.
For the minimization problem to make sense as an {\em approximation},
we need precisely this condition.
In \cite{GEORGIOU_LINDQUIST_KULLBACKLEIBLER} it is shown that, when $A$ is singular,
the zeroth-order moment of all the spectra $\Phi$ compatible
with the constraint $\int G \Phi G^\ast = \Sigma$ is constant,
say $\int\Phi \equiv \alpha$.
Without either of the singularity of $A$ or the equality of the zeroth-order
moments of $\Psi$ {\em and} all of the $\Phi$'s, 
it is not clear at all if $\KLD$ serves as a pseudo-distance,
even if the minimization problem continues to be valid.
That is, it is not clear whether we can speak about ``approximation'' anymore.
In view of this consideration, we require from now on
that $A$ has at least one eigenvalue at the origin, and
that $\Psi$ is rescaled accordingly in order to obtain $\int\Psi=\alpha$.
(Hence, what we approximate is the ``shape'' of $\Psi$, not $\Psi$ itself.)\\
If $A$ is non-singular, it is still possible to consider a weighted version of the Kullback-Leibler pseudo-distance in such a way that the problem maintain the meaning of spectral approximation.
\end{rem}

To simplify the writing, we can, without loss of generality, normalize $\Sigma$ and $\Psi$.
Indeed, if $\Sigma\neq I$, it suffices to replace $G$ by
$G':=\Sigma^{-1/2}G$ and $(A,B)$ with
$(A'=\Sigma^{-1/2}A\Sigma^{1/2},B'=\Sigma^{-1/2}B)$ to obtain an equivalent problem
where $\Sigma= I$.
%
%

In a similar fashion, if $\int\Phi \equiv \int\Psi = \alpha \neq 1$
(compare with Remark \ref{REMARK_KL_POSITIVITY}), let $\Psi':=\Psi/\alpha$
and $G' = \sqrt{\alpha}\ G$. Then, to any solution $\Phi$ to the
moment problem $\int G \Phi G^\ast = \Sigma$ there corresponds a solution
$\Phi'$ to the problem $\int G' \Phi' G'^\ast = \Sigma$,
where $\Phi' = \Phi/\alpha$.
It is immediate to check that
$$
\KLD(\Psi'\|\Phi') = \KLD(\Psi\|\Phi) / \alpha,
$$
which ensures that the positivity of the pseudo-distance is preserved.
Therefore, we can assume that $\int\Psi = 1$.

The first issue one needs to worry about is
existence of $\Phi\in C_+(\T)$ satisfying  constraint (\ref{2}).
It has been shown that the following conditions are  equivalent \cite{GEORGIOU_LINDQUIST_KULLBACKLEIBLER}:
\begin{enumerate}
\item
The family of $\Phi$ satisfying  constraint (\ref{2}) is nonempty.
\item
 there exists
$H\in
\C^{1\times n}$ such that
\beq
I-A  A^*=BH+H^*B^*,
\eeq
\item
the following rank condition holds
\begin{equation}\label{3}{\rm rank}\left(\begin{array}{cc}I-A
A^*&B\\B^*&0\end{array}\right)={\rm rank}\left(\begin{array}{cc}0&B\\B^*&0
\end{array}\right)
\end{equation}
\end{enumerate}
A fourth equivalent condition is based on a linear operator that will play a crucial role in  the rest of the paper: Let $\Hermitian_n$  be the space of Hermitian matrices of dimension $n$, and consider the linear operator
\beq
\begin{array}{cccc}
\Gamma: & C(\T)& \longrightarrow & \Hermitian_n\\
&\Phi & \mapsto & \int G\Phi  G^*
\end{array}
\eeq
It is clear (recall that $\Sigma =I$) that there exists $\Phi\in C(\T)$ satisfying (\ref{2}) if and only if 
\beq\label{iinrga}
I\in \Rgamma.
\eeq
Indeed, it has been shown \cite{FERRANTE_PAVON_RAMPONI_HELLINGERVSKULLBACK} that condition (\ref{iinrga}) is necessary and sufficient for the family of $\Phi$ in (\ref{2}) to be nonempty.
Thus, condition (\ref{iinrga}) will be a standing assumption for this paper.
We endow  the space $\Hermitian_n$  of Hermitian matrices with the inner product
$$
\left<P,Q\right>:=\tr(PQ).
$$
The orthogonal complement of 
$\Rgamma$ (defined with respect to this inner product)  has been shown in  \cite{FERRANTE_PAVON_RAMPONI_HELLINGERVSKULLBACK} to be given by
\bea\nn
&\Rgamma^\perp=\left\{X\in\Hermitian_n:\right.\hspace{4cm}\\
&\hspace{2cm}\hfill\left. \ G^\ast(\ejth)X G(\ejth) = 0, \ \forall \vartheta \in [0, 2\pi]\right\}.
\eea
Notice that, in view of (\ref{iinrga}), we have
$$\tr(X)=\left<X,I\right>=0,\quad \forall X\in\Rgamma^\perp.$$
\subsection{Variational analysis}
To solve Problem \ref{AP}, we consider a {\em matrix Lagrange multiplier}
$\Lambda\in\Hermitian_n$ satisfying  $G^*\Lambda G>0$  on all of $\T$, and
define the {\em Lagrangian  functional}
\bea\nn
L(\Phi,\Lambda)&:=&\KLD(\Psi\|\Phi)+\left<\Lambda,
\int G\Phi G^*-I\right>\\
&=&\KLD(\Psi\|\Phi)+\int  G^*\Lambda
G\Phi-\tr(\Lambda).\label{4}
\eea
This functional is easily seen to be strictly convex.
Therefore {\em unconstrained} minimization of 
$L(\Phi,\Lambda)$ can be achieved by annihilating the directional derivative
$\DD(L(\Phi,\Lambda),\delta\Phi)$ along all directions $\delta\Phi\in C(\T)$.
In this way we get the following form for the optimal solution as a function of the Lagrange multiplier $\Lambda$.
\beq\label{phioptfdl}
\Phi_\Lambda=\frac{\Psi}{G^* \Lambda G}
\eeq
Now, it is clear that if  $\Lambdaopt=\Lambdaopt^*$ satisfies 
\bea\label{condlam}
\stoa\label{9} &&G^*\Lambdaopt G>0, \quad\forall
\ejth\in\T,\\ 
\stob\label{10b}
&&\int G\frac{\Psi}{G^*\Lambdaopt G}G^*=I\label{10}
\eea\rst
then 
\beq\label{phiopt}\Phi_\circ:=\Phi_{\Lambdaopt}=\frac{\Psi}{G^* \Lambdaopt G}
\eeq
is  optimal for Problem \ref{AP}.
As for many optimization problems, the most delicate issue is {\em existence}. This issue has been addressed in \cite{GEORGIOU_LINDQUIST_KULLBACKLEIBLER,FERRANTE_PAVON_RAMPONI_FURTHERRESULTS} where the following result has been proven.
\begin{theo} \label{Opt} There exist matrices 
$\Lambdaopt=\Lambdaopt^*$  such that (\ref{condlam}) hold.
For any such a $\Lambdaopt$,  $\Phi_\circ$ given by (\ref{phiopt}) is the unique solution of the Approximation Problem
(\ref{AP}).
\end{theo}
\begin{rem}\label{remrgort}
Notice that, since Problem (\ref{AP}) admits a unique solution, 
if $\Lambdaopt$ and $\Lambdaopt'$ are two matrices satisfying conditions (\ref{condlam}), then $\frac{\Psi}{G^* \Lambdaopt G}=\frac{\Psi}{G^* \Lambdaopt' G}$ so that we clearly have
$
G^* (\Lambdaopt-\Lambdaopt') G\equiv 0$, or equivalently,
 $\Lambdaopt-\Lambdaopt'\in\Rgamma^\perp$.
Conversely, it is clear that if $\Lambdaopt$  satisfies conditions (\ref{condlam})
then $(\Lambdaopt+X)$ also satisfies conditions (\ref{condlam}) for any $X\in\Rgamma^\perp$.
Thus, the family ${\cal L}_{\circ}$ of all  solutions of (\ref{condlam}) is an affine space that may be parametrized in terms of an arbitrary solution $\Lambdaopt$, as 
\beq\label{callz}
{\cal L}_{\circ}=\{\Lambdaopt+X;\ X\in(\Rgamma)^\perp\}.
\eeq
Moreover, for any $\Lambdaopt$ satisfying (\ref{condlam}), we have
$\tr\int G\frac{\Psi}{G^*\Lambdaopt G}G^*\Lambdaopt =\tr\Lambdaopt $
and, using the cyclic property of the trace we immediately get:
$$\tr\Lambdaopt=1.$$ 
\end{rem}
By duality theory, a $\Lambdaopt$ satisfying (\ref{condlam}) may be computed by 
maximization of the dual functional
$$\Lambda\mapsto \inf\{L(\Phi,\Lambda) | \Phi\in C_+(\T)\}.
$$

The latter, in view of the previous discussion, may be explicitly written as:
\begin{equation}\label{dual}\Lambda\mapsto L\left(\frac{\Psi}{G^*\Lambda
G},\Lambda\right)=\int\Psi\log G^*\Lambda G -\tr (\Lambda) +\int\Psi.
\end{equation}

Consider now the maximization of the dual functional
(\ref{dual}) over the set
\begin{equation}\label{lam} {\cal L}_+:=\{\Lambda=\Lambda^* |G^*\Lambda
G>0,
\forall \ejth\in\T\}.
\end{equation}
 Let$$\J_\Psi(\Lambda):=-\int\Psi\log G^*\Lambda G +\tr (\Lambda).
$$
The dual problem is then equivalent to
\begin{equation} \label{DP}{\rm minimize}\quad\{\J_\Psi(\Lambda)
|\Lambda\in {\cal L}_+\}.
\end{equation}
 
As discussed in the Introduction, the bottleneck of the whole theory and of its numerous applications is now the numerical computation of a  $\Lambdaopt$ satisfying  (\ref{condlam}) or, equivalently, solving (\ref{DP}). To this aim the following algorithm has been proposed in \cite{PAVON_FERRANTE_ONGEORGIOULINDQUIST} and further discussed in \cite{FERRANTE_PAVON_RAMPONI_FURTHERRESULTS}.

\subsection{Iterative algorithm}\label{iteralg}
For $\Lambda\geq0$, let
\beq\label{THETA}
\Theta(\Lambda):=\int \Lambda^{1/2}G\left[\frac{\Psi}{G^*\Lambda G}\right]G^*\Lambda^{1/2}.
\eeq
It has been shown in \cite{PAVON_FERRANTE_ONGEORGIOULINDQUIST} that $\Theta$ is a map from {\em density   matrices} to density    matrices, i.e. if $\Lambda$ is a positive semi-definite Hermitian matrix with trace equal to $1$, then  
$\Theta(\Lambda)$ has the same properties.   Density matrices have long  been studied in statistical quantum  mechanics, representing quantum states in the presence of uncertainty \cite{VONNEUMANN_MATHEMATICALFOUNDATIONS,NIELSEN_CHUANG_QUANTUMCOMPUTATION}. 
Moreover, $\Theta$ maintains positive definiteness, i.e., if $\Lambda>0$, then $\Theta(\Lambda)>0.$
 In addition to this, the following holds:  
\begin{prop}
The matrix $\Theta(\Lambda)$ has the same rank and, indeed, the same kernel of the matrix $\Lambda$.
\end{prop}
\IEEEproof
By taking into account that $\ker (\Lambda^{1/2})=\ker (\Lambda)$, it is immediate to check that if
$v\in\ker\left[\Lambda\right]$ then $v\in\ker\left[\Theta(\Lambda)\right]$.
Conversely, it is sufficient to prove that 
\beq\label{fmn}
\int G\left[\frac{\Psi}{G^*\Lambda G}\right]G^*>0.
\eeq
Indeed, if  this is the case and  $v\in\ker\left[\Theta(\Lambda)\right]$ then
$\Lambda^{1/2}v=0,$ so that $v\in\ker\left[\Lambda\right]$.
To prove (\ref{fmn}), we observe that $\frac{\Psi}{G^*\Lambda G}$ is continuous and strictly positive on $\T$ and hence has a positive minimum there.
It is therefore sufficient to show that $\int GG^*$ is positive definite. The latter integral is the steady-state state covariance of the filter $G$ driven by normalized white noise, i.e., the unique solution $\Xi$ of the discrete-time Lyapunov equation 
$\Xi -A\Xi A^* =BB^*$. In view of controllability of the pair $(A,B)$, it is clear that 
$\Xi$ is positive definite.
\cvd

Consider the sequence $\{\Lambda_{k}\}$ produced by the following iteration
\beq\label{algori} 
\Lambda_{k+1}=\Theta(\Lambda_{k}),
\eeq
with an arbitrary initial condition $\Lambda_{0}>0$.
Notice that,   since  each $\Lambda_{k}$ produced by the previous algorithm is positive definite and has trace equal to $1,$ we also have
$\Lambda_{k}\leq I$ $\forall k>0$.   If the  sequence $\{\Lambda_{k}\}$ converges to a limit point $\hat{\Lambda}>0$ then such a  $\hat{\Lambda}$ 
is a fixed point for the map $\Theta$ in (\ref{THETA}):
\beq\label{THETAFP}
\hat{\Lambda}:=\int \hat{\Lambda}^{1/2}G\left[\frac{\Psi}{G^*\hat{\Lambda} G}\right]G^*\hat{\Lambda}^{1/2}.
\eeq
By multiplying the latter by $\hat{\Lambda}^{-1/2}$ on both sides, it is clear that 
$\hat{\Lambda}$ satisfies (\ref{condlam})  and hence provides a solution of Problem \ref{AP}.

Notice that, even if all $\Lambda_k$ are positive definite, it may happen that   the  sequence $\{\Lambda_{k}\}$ converges to a limit point $\hat{\Lambda}_s$ which is positive semidefinite but singular. In this case, it is not guaranteed that   $\hat{\Lambda}_s$ satisfy (\ref{condlam}).

We observe that if $\Lambdaopt>0$ is a fixed point of $\Theta$, then
for any $\Lambda_{\perp}\in(\Rgamma)^\perp$,
$\Lambdaopt+\Lambda_{\perp}$ is also a fixed point of $\Theta$, as long as
$\Lambdaopt+\Lambda_{\perp}\geq 0$.
In fact, in view of Remark \ref{remrgort}, we have
\beq\label{serveperprop}
\Theta(\Lambdaopt+\Lambda_{\perp})=(\Lambdaopt+\Lambda_{\perp})^{1/2}I(\Lambdaopt+\Lambda_{\perp})^{1/2}=\Lambdaopt+\Lambda_{\perp}.
\eeq

In a wide series of simulations, we have observed that $\Lambda_{k}$ always converges to a limit point. In only one case such a limit point was a singular matrix. Also in that case, however, the limit point   satisfies (\ref{condlam}) and hence  provides a solution of Problem \ref{AP}.

\section{Proof of convergence}\label{sec:POC}

  In this section we prove the main contribution of the paper, namely
that the intersection between the affine family of  solutions of (\ref{condlam}) and the cone of positive definite matrices is a locally asymptotically stable manifold for the iteration (\ref{THETA}).

\subsection{Existence of a positive definite $\Lambdaopt$}

Once again, the first issue that must be addressed is an existence result: We have to show that
\beq\label{callzp}
{\cal L}_{\circ+}:=\{\Lambdaopt=\Lambdaopt^*>0:\ \Lambdaopt\in{\cal L}_{\circ}\}\neq\emptyset,
\eeq
where ${\cal L}_{\circ}$ is defined in (\ref{callz}).
To this aim we need a preliminary result in the same vein of Lemma 9 in \cite{GEORGIOU_LINDQUIST_KULLBACKLEIBLER}. The latter has been established in a slightly different setting and using an abstract functional-analytic approach.
We will instead use a direct algebraic approach that provides a constructive proof.
\begin{lem}\label{analoglemma9}
If $G^*\Lambdaopt G>0, \ \forall\,\ejth\in\T$, then there exists
a vector
$C_\circ\in\C^{n\times  1}$ such that
$G^*\Lambdaopt G=G^*C_\circ C_\circ^* G$.
\end{lem}
\IEEEproof
As shown in Lemma \ref{lemapp} in the Appendix,   we can obtain a decomposition 
\beq
G^*\Lambda G=W^*W,
\eeq
where   the (right) spectral factor $W(z)$ is given by (\ref{wlem}).
  Denoting by $P$ the stabilizing solution of the Riccati equation (\ref{AREL}), by using (\ref{last}), $W$ may be explicitly expressed in the form 
\beq\label{wrew}
W=(B^{\ast}PB)^{-1/2}B^{\ast}P\left(A(zI-A)^{-1}+I\right)B.
\eeq
It is immediate to check that $A(zI-A)^{-1}+I=z(zI-A)^{-1}$
so that
\beq\label{wrerew}
W=z(B^{\ast}PB)^{-1/2}B^{\ast}P(zI-A)^{-1}B.
\eeq
and thus
\beq
G^*\Lambda G=W^*W=W_1^*W_1,
\eeq
with 
\beq
W_1:=z^{-1}W=(B^{\ast}PB)^{-1/2}B^{\ast}P(zI-A)^{-1}B.
\eeq
Therefore, the vector $C_\circ$ indeed exists and may be explicitly expressed as
$C_\circ=\left((B^{\ast}PB)^{-1/2}B^{\ast}P\right)^*$.
\cvd
\begin{theo}
The set ${\cal L}_{\circ+}$ defined by (\ref{callzp}) is nonempty and it is an open convex subset of the affine space ${\cal L}_{\circ}$.
\end{theo}
\IEEEproof
Let $\Lambdaopt\in{\cal L}_{\circ}$ (recall that Theorem \ref{Opt} guarantees that ${\cal L}_{\circ}\neq\emptyset$) so that 
$G^*\Lambdaopt G>0, \ \forall\,\ejth\in\T$.
From Lemma \ref{analoglemma9} we know that this implies the existence of a vector
$C_\circ\in\C^{n\times  1}$ such that
$G^*\Lambdaopt G=G^*C_\circ C_\circ^* G$.\\
On the unit circle $\T$, $G^*\Lambdaopt G$ is continuous and positive and 
hence
$$
\mu:=\min\{G(z)^*\Lambdaopt G(z):\ z\in\T\}>0.
$$
Similarly, on the unit circle $\T$, $G^* G$ is continuous and 
hence
$$
\nu :=\max\{G(z)^*G(z):\ z\in\T\}
$$
is finite.
Let $\varepsilon:=\frac{\mu}{4\nu }$.
Clearly, $\forall\  z\in\T$,
\bea
\nn
G(z)^*(\frac{1}{2}\Lambdaopt-\varepsilon I)G(z)&=&\frac{1}{2}G(z)^*\Lambdaopt G(z)
-\varepsilon G(z)^* G(z)\\
&\geq &\frac{\mu}{2}-\frac{\mu}{4\nu }\nu =
\frac{\mu}{4}>0.
\eea
Hence, exploiting again Lemma \ref{analoglemma9}, we conclude that there
exists  $C_1\in\C^{n\times  1}$ such that
$$G^*(\frac{1}{2}\Lambdaopt-\varepsilon I) G=G^*C_1 C_1^* G.$$
Therefore we have
\bea
\nn
G^*\Lambdaopt G\!\!\!&=\!\!\!&\frac{1}{2}G^*C_\circ C_\circ^* G+\frac{1}{2}G^*C_\circ C_\circ^* G +\varepsilon G ^* G -\varepsilon G ^* G \\
\nn
&=\!\!\!&
G^*(\frac{1}{2}C_\circ C_\circ^*+\varepsilon I)G +\frac{1}{2}G^*C_\circ C_\circ^* G -\varepsilon G ^* G \\
\nn
&=\!\!\!&
G^*(\frac{1}{2}C_\circ C_\circ^*+\varepsilon I)G +\frac{1}{2}G^*\Lambdaopt G -\varepsilon G ^* G \\
\nn
&=\!\!\!&
G^*(\frac{1}{2}C_\circ C_\circ^*+\varepsilon I)G +G^*(\frac{1}{2}\Lambdaopt-\varepsilon I) G\\
\nn
&=\!\!\!&
G^*(\frac{1}{2}C_\circ C_\circ^*+\varepsilon I)G +G^*C_1 C_1^* G\\
\nn
&=\!\!\!&
G^*(\frac{1}{2}C_\circ C_\circ^*+\varepsilon I+C_1 C_1^*)G=G^*\Lambda_{\circ+}G,
\eea
where $\Lambda_{\circ+}:=\frac{1}{2}C_\circ C_\circ^*+\varepsilon I+C_1 C_1^*$ is clearly positive definite and hence $\Lambda_{\circ+}\in{\cal L}_{\circ+}$.
The fact that ${\cal L}_{\circ+}$   is an open convex subset of  ${\cal L}_{\circ}$ is an immediate consequence of the fact that the cone of positive definite matrices is open and convex together with the fact that
${\cal L}_{\circ}$ is an affine space.
\cvd

\subsection{Linearization}

  Given that ${\cal L}_{\circ+}$ is non-empty,  we can now pick a point $\Lambdaopt\in{\cal L}_{\circ+}$ and analyze the map $\Theta$ in a neighborhood of $\Lambdaopt$.
To this aim we linearize the map $\Theta,$ namely we compute the directional derivative of $\Theta$ at $\Lambdaopt$  in the direction specified by an arbitrary Hermitian matrix $X$:
$$
D(\Theta(\Lambdaopt),X):=\lim_{\varepsilon\rightarrow0}\frac{\Theta(\Lambdaopt+\varepsilon X)-\Theta(\Lambdaopt)}{\varepsilon}.
$$
In order to find an explicit form for this derivative, we first need an expression for 
$D(\Lambdaopt^{1/2},X)$.
\subsubsection{Derivative of the matrix square root}
For a given function $\varrho:\Hermitian_n\rightarrow \Hermitian_n$, let us take the directional derivative of $\left(\varrho(\Lambda)\right)^2$ in the direction $X$.   The chain rule gives:
\begin{displaymath}
D(\varrho(\Lambda)^2; X) = D(\varrho(\Lambda);X)\ \varrho(\Lambda) + \varrho(\Lambda)\ D(\varrho(\Lambda);X)
\end{displaymath}
Now if $\varrho(\Lambda) = \Lambda^{1/2}$, we have
$\varrho(\Lambda)^2=\Lambda$ so that clearly $D(\varrho(\Lambda)^2; X)=D(\Lambda; X)=X$. 
In conclusion, we get that the derivative $\DD(\Lambda^{1/2};X)$
is the solution of the following Lyapunov equation:
\begin{equation}
\label{SQDERIV}
\DD(\Lambda^{1/2};X)\ \Lambda^{1/2} + \Lambda^{1/2}\ \DD(\Lambda^{1/2};X) = X
\end{equation}
%

\subsubsection{Derivative of $\Theta$}


\newcommand{\blambda}{\bar{\Lambda}}

Let us take the variation of (\ref{THETA}) in a direction $X$.
  By applying the chain rule we get:
\begin{equation}
\label{VARTHETA1}
\begin{split}
\DD(\Theta(\Lambda);X) 
=\  & \DD(\Lambda^{1/2};X) \int \frac{G\Psi G^\ast}{G^\ast \Lambda G} \Lambda^{1/2} \\
& + \Lambda^{1/2} \DD\left(\int \frac{G\Psi G^\ast}{G^\ast \Lambda G} ; X\right) \Lambda^{1/2} \\
&+ \Lambda^{1/2} \int \frac{G\Psi G^\ast}{G^\ast \Lambda G} \ \DD(\Lambda^{1/2};X) \\
= \ & \DD(\Lambda^{1/2};X) \int \frac{G\Psi G^\ast}{G^\ast \Lambda G} \Lambda^{1/2} \\
& - \Lambda^{1/2} \int \frac{G\Psi G^\ast}{G^\ast \Lambda G} \frac{G^\ast X G}{G^\ast \Lambda G} \Lambda^{1/2} \\
&+ \Lambda^{1/2} \int \frac{G\Psi G^\ast}{G^\ast \Lambda G} \ \DD(\Lambda^{1/2};X) 
\end{split}
\end{equation}
We now compute the latter expression at $\Lambda = \Lambdaopt
$ and take (\ref{10b}) into account. This yields:
\begin{equation}
\label{VARTHETA2}
\begin{split}
\DD(\Theta(\Lambdaopt
);X) 
=\ & \DD(\Lambdaopt
^{1/2};X) I\ \Lambdaopt
^{1/2} + \Lambdaopt
^{1/2} \ I\ \DD(\Lambdaopt
^{1/2};X)
 \\
& - \Lambdaopt
^{1/2} \int \frac{G\Psi G^\ast}{G^\ast \Lambdaopt
 G} \frac{G^\ast X G}{G^\ast \Lambdaopt
 G} \Lambdaopt
^{1/2},
\end{split} 
\end{equation}
which, by property (\ref{SQDERIV}), may be rewritten as
\begin{equation}
\label{VARTHETA}
\DD(\Theta(\Lambdaopt
);X) 
= X - \Lambdaopt
^{1/2} \int \frac{G\Psi G^\ast}{G^\ast \Lambdaopt
 G} \frac{G^\ast X G}{G^\ast \Lambdaopt
 G} \Lambdaopt
^{1/2}.
\end{equation}
Let us define the linear map ${\mathscr M}:\Hermitian_n\longrightarrow\Hermitian_n $ as the above derivative:
\begin{equation}
\label{M}
 {\mathscr M}(X) :=  
X - \Lambdaopt
^{1/2} \int \frac{G\Psi G^\ast}{G^\ast \Lambdaopt
 G} \frac{G^\ast X G}{G^\ast \Lambdaopt
 G}\ \Lambdaopt
^{1/2} 
\end{equation}
The linear map ${\mathscr M}$ is then the derivative of $\Theta$ computed at a given fixed point $\Lambdaopt$.

Adopting a system-theoretic approach, we can consider the sequence of increments
$\{X_k\}$ with $X_k=\Lambda_k-\Lambdaopt$ and the {\em linear} system
$X_{k+1} =  {\mathscr M}(X_k)$ as the linear approximation
of the nonlinear discrete-time system $\Lambda_{k+1} = \Theta(\Lambda_k)$ in the neighborhood of
its equilibrium point $\Lambdaopt$.
If all the eigenvalues of $ {\mathscr M}$ lied in the open unit circle, then we could immediately conclude that $\Lambdaopt$ is asymptotically stable.
However, this is {\em not} the case. In fact, it is immediate to check that 
\beq\label{miden}
{\mathscr M}(X_\perp)=X_\perp,\quad \forall\ X_\perp\in\Rgamma^\perp
\eeq  so that $ {\mathscr M}$ restricted to $\Rgamma^\perp$ is the identity operator.
We thus need a more sophisticated analysis.

\subsection{Properties and spectrum of $ {\mathscr M}$}

First, notice that $ {\mathscr M}$ maps $\Hermitian_n$ in itself but it is not self-adjoint
  (with respect to the inner product  defined in $\Hermitian_n$ by $\left< X,Y \right> = \tr X Y$).
Indeed, given $X,Y\in\Hermitian_n$, it may happen that
\begin{displaymath}
\left<  {\mathscr M}(X),Y \right> \neq \left< X,  {\mathscr M}(Y) \right> 
\end{displaymath}
so it is not {\em a priori} true that 
 the eigenvalues of $ {\mathscr M}$ are real and that 
the eigenmatrices of $ {\mathscr M}$ span
the whole space 
$\Hermitian_n$.\\
A second observation is stated in the following result.
\begin{lem}
For any $X\in\Hermitian_n$,
$\tr  {\mathscr M}(X)=0$.
\end{lem}
\IEEEproof We have
\begin{equation}
\label{TRACEM}
\begin{split}
\tr  {\mathscr M}(X) 
&= \tr X - \int \Psi \frac{G^\ast \Lambdaopt
 G}{G^\ast \Lambdaopt
 G} \frac{G^\ast X G}{G^\ast \Lambdaopt
 G}\\ 
& = \tr X - \int \Psi \frac{G^\ast X G}{G^\ast \Lambdaopt
 G} \\
&= \tr X - \tr \int \frac{G \Psi G^\ast}{G^\ast \Lambdaopt
 G}\ X 
 = \tr X - \tr I\ X 
 = 0
\end{split}
\end{equation}
\cvd

We are now ready to analyze the spectrum of $ {\mathscr M}$. Let $Y$ be an eigenmatrix of $ {\mathscr M}$ and  $\alpha$ be the corresponding eigenvalue, namely $Y$ is a non-zero Hermitian matrix such that  $ {\mathscr M}(Y) = \alpha\ Y$.
Due to (\ref{TRACEM}), we have
\begin{equation}
\label{EIGENZERO}
\alpha\ \tr Y = 0
\end{equation}
Thus, we have the following corollary.
\begin{cor}\label{cortace}
Let $Y$ be any eigenmatrix of ${\mathscr M}$ and assume $\tr Y\neq 0.$ Then, the corresponding eigenvalue is zero.
\end{cor}

Notice that $\Lambdaopt$ is one such eigenmatrix. Indeed:
\begin{displaymath}
\begin{split}
 {\mathscr M}(\Lambdaopt
) &=  
\Lambdaopt
 - \Lambdaopt
^{1/2} \int \frac{G\Psi G^\ast}{G^\ast \Lambdaopt G} 
\frac{G^\ast \Lambdaopt G}{G^\ast \Lambdaopt G}\ 
\Lambdaopt^{1/2} \\
&= \Lambdaopt - \Lambdaopt^{1/2} 
\int \frac{G\Psi G^\ast}{G^\ast \Lambdaopt G}\ 
\Lambdaopt^{1/2} \\
&= \Lambdaopt
 - \Lambdaopt
^{1/2}\ I\ \Lambdaopt^{1/2}=0.
\end{split}
\end{displaymath}

Let $Y$ be an eigenmatrix of $ {\mathscr M}$ and 
$\alpha$ be the corresponding eigenvalue. We want to compute bounds for $\alpha$.
In view of Corollary \ref{cortace}, we can
assume $\tr Y = 0$.
We have:
$$
\alpha Y  =  Y - \Lambdaopt
^{1/2} \int \frac{G\Psi G^\ast}{G^\ast \Lambdaopt
 G} \frac{G^\ast Y G}{G^\ast \Lambdaopt
 G}\ \Lambdaopt
^{1/2} $$
or, equivalently,
$$
(1 - \alpha)Y   =  \Lambdaopt
^{1/2} \int \frac{G\Psi G^\ast}{G^\ast \Lambdaopt
 G} \frac{G^\ast Y G}{G^\ast \Lambdaopt
 G}\ \Lambdaopt
^{1/2}.
$$
Since $\Lambdaopt>0$, we can multiply both members  by  
$\Lambdaopt^{-1/2}$ on the left side, and by $\Lambdaopt^{-1/2}Y$ on the right side. This yields
$$
(1 - \alpha)\Lambdaopt^{-1/2}Y\Lambdaopt^{-1/2}Y  =  
\int \frac{G\Psi G^\ast}{G^\ast \Lambdaopt G} 
\frac{G^\ast Y G}{G^\ast \Lambdaopt G}Y,$$
which, by taking the trace on both members and exploiting the cyclic property of the trace, implies
$$
(1 - \alpha)\tr\left[\Lambdaopt^{-1/2}Y\Lambdaopt^{-1/2}Y\right]  =  \int \Psi \frac{(G^\ast Y G)^2}{(G^\ast \Lambdaopt
 G)^2}. $$
We now observe that 
$$\tr\left[\Lambdaopt^{-1/2}Y\Lambdaopt^{-1/2}Y\right]
=\tr \left[(\Lambdaopt ^{-1/4}Y\Lambdaopt^{-1/4})^2\right]$$
which is strictly positive because it is the square of the Frobenius norm
of the nonzero matrix $\Lambdaopt ^{-1/4}Y\Lambdaopt^{-1/4}$.
In conclusion we get
\beq\label{alppos}
(1 - \alpha) =  \frac{ \int \Psi \frac{(G^\ast Y G)^2}{(G^\ast \Lambdaopt
 G)^2} }{ \tr \left[(\Lambdaopt ^{-1/4}Y\Lambdaopt
^{-1/4})^2\right]  } .
\eeq
The right-hand side of (\ref{alppos}) is clearly real and non-negative.
Indeed, it vanishes if and only if $G(z)^\ast Y G(z)$ is identically zero on $\T$ or, equivalently if an only if  $Y \in (\Range\Gamma)^\perp$.
The following theorem is thus proven.
\begin{theo}\label{theo32}
All the eigenvalues of the map $ {\mathscr M}$ are real. For any eigenmatrix of $ {\mathscr M}$ that is {\em not} in
$ (\Range\Gamma)^\perp$, the corresponding eigenvalue is strictly smaller than $1$. On the space $(\Range\Gamma)^\perp$, $ {\mathscr M}$ acts as the identity operator.
\end{theo}
\begin{rem}\label{remvectoriz}
The above theorem may be interpreted as follows.
We define $\bar{x}:={\rm vec}(X)$ as the column vector (with $n^2$ entries) obtained by stacking the columns of $X$ one over the other.
Let $V$ be a matrix whose columns form a basis for $\{x={\rm vec}(X):\  X\in(\Range\Gamma)^\perp\}$. Let $W$ be such that 
\beq\label{wt}
T:=[V\mid W]
\eeq is nonsingular.
Let $x:=T^{-1}\bar{x}=T^{-1}{\rm vec}(X)$. Clearly, $x$ is a coordinate representation of $X$.
Theorem \ref{theo32} states that, with respect to these coordinates
the linear map ${\mathscr M}$ is represented by a matrix $M$
of dimension $n^2\times n^2$ with the following structure: 
$
M=\bmat{cc} I_{n_\perp} & M_{12}\\0& B\emat,
$
where $n_\perp$ is the dimension of $(\Range\Gamma)^\perp$ and 
\beq\label{sigb}
\sigma(B)\subset  (-\infty, 1).
\eeq
Clearly, since $I_{n_\perp}$ and 
$B$ have disjoint spectra, we can select $W$ in
(\ref{wt}) such that $M_{12}=0$, i.e. $M$ ha the structure
\beq\label{matrep}
M=\bmat{cc} I_{n_\perp} & 0\\0& B\emat.
\eeq
It remains to establish a lower bound for the spectrum of $B$.
\end{rem}

\subsection{Eigenvalues of $ {\mathscr M}$ are non-negative}

In order to provide a lower bound for the spectrum, we shall consider the linearized map as the generator of a continuous-time semigroup evolution, for which a key spectral property will be derived.
\subsubsection{$ {\mathscr M}$ is the opposite of a {\em Lindblad} generator}
We observe that the  operator $ {\mathscr M}$ may be written in the form
\bea
\label{Mlim}
\nn
 {\mathscr M}(X) &\!\!\!\!=\!\!\!\!&  
X \! - \! 
\int \left[\Lambdaopt^{1/2}G\frac{\Psi^{1/2}}{G^\ast \Lambdaopt G}G^\ast\right]X
\left[G\frac{\Psi^{1/2}}{G^\ast \Lambdaopt G}G^\ast\Lambdaopt^{1/2}\right]\\
&\!\!\!\!=\!\!\!\!&
X-\int LXL^\ast
\eea
where $L:=\Lambdaopt^{1/2}G\frac{\Psi^{1/2}}{G^\ast \Lambdaopt G}G^\ast$.
It is immediate to check that $L^\ast L=G\frac{\Psi}{G^\ast \Lambdaopt G}G^\ast$,
so that
clearly $\int L^\ast L =I$.
Therefore, we can write write $- {\mathscr M}(X)$ as a (generalized) {\em Lindblad}   generator \cite{lindblad}:\footnote{It is remarkable to notice that, in the framework of quantum statistical mechanics, it has been shown by Lindblad \cite{lindblad} that any trace-preserving, strongly-continuous semigroup of completely positive maps from density operators to density operators has a generator which can be written as the sum of an Hamiltonian (Liouvillian) term and number of terms of the form of the integrand in \eqref{lindbladterm}. Such Markov semigroups have   long been studied for their relevance  to many aspects of quantum theory and thermodynamics \cite{alicki-lendi,petruccione}. In this setting, their spectral properties have been investigated from an operator-theoretic standpoint. In order to avoid to overburden this paper with an unnecessary and rather technical detour, we choose here to prove the needed results by means of a linear algebraic tools.}
\beq\label{lindbladterm}
- {\mathscr M}(X)= \int LXL^\ast -\frac{1}{2}[XL^\ast L+L^\ast LX]
\eeq

\subsubsection{A continuous-time evolution}
We now consider the following continuous-time linear system 
\beq\label{ctev}
\dot{X}(t)=- {\mathscr M}(X(t))=\int LX(t)L^\ast -\frac{1}{2}[X(t)L^\ast L+L^\ast LX(t)].
\eeq
with state space being the set of traceless Hermitian matrices (notice that, since 
$ \tr{\mathscr M}(X)=0$, evolution (\ref{ctev}) is trace preserving). 
This will be helpful in proving the following 
\begin{theo}\label{theo33}
All the eigenvalues of the map $ {\mathscr M}$ are non-negative.
\end{theo}
\IEEEproof
Let $\alpha$ be an eigenvalue of $ {\mathscr M}$ and $Y$ be the corresponding eigenmatrix, so that
the state trajectory generated by system (\ref{ctev}) with  initial condition
$X(0)=Y$ is $X(t)={\rm e}^{-\alpha t}Y$.
We denote  by $\|Y\|_1$ the sum of the absolute values of the eigenvalues of $Y$, i.e. 
$\|Y\|_1:={\displaystyle \sum_{\lambda\in\sigma(Y)}}|\lambda|$.

Let $Y_P\geq 0$ and $Y_N \geq 0$ be the {\em positive and negative parts} of $Y$ defined as follows: Let $T^\ast YT=D=D_P-D_N$, where $T^\ast=T^{-1}$, $D$ is a diagonal matrix, $D_P$ is obtained from $D$ by annihilating the negative entries and  $D_N:=-(D-D_P)$.
Define $Y_P:=TD_PT^\ast$ and $Y_N:=TD_NT^\ast$. Define also the orthogonal projection
 $\Pi_P:=TO_PT^\ast$ ($\Pi_N:=TO_NT^\ast$), where $O_P$ ($O_N$) is the matrix obtained from
 $D_P$ ($D_N$) by setting to $1$ all the non-zero entries.
Clearly,
\beq\label{pd1}
\|Y \|_1=\tr[Y_P + Y_N ],\quad Y =Y_P -Y_N .
\eeq
Moreover,
\beq\label{pd2}
Y_P =\Pi_P Y ,\quad Y_N =-\Pi_N Y , \quad Y_P  Y_N = Y_N  Y_P =0.
\eeq
Recall now that, in view of Corollary \ref{cortace}, we can assume $\tr[Y]= 0$. Thus, taking into account   (\ref{pd1}) and (\ref{pd2}), we have
$0=\tr[Y]=\tr[(\Pi_P+ \Pi_N)Y]$ so that
$0=\tr[X(t)]=\tr[(\Pi_P+ \Pi_N)X(t)]$.
Hence 
$$ \tr[\Pi_P\frac{d}{dt} X(t)]=- \tr[\Pi_N\frac{d}{dt} X(t)].$$
It is now easy to see that
\bea\nn\frac{d}{dt} \|X(t)\|_1
&\!\!\!=\!\!\!&
\frac{d}{dt} \|{\rm e}^{-\alpha t}Y\|_1=\frac{d}{dt} \tr[{\rm e}^{-\alpha t}(Y_P+Y_N)]\\
\nn
&\!\!\!=\!\!\!& \tr[\frac{d}{dt} ({\rm e}^{-\alpha t}(\Pi_P-\Pi_N)Y)]\\
\nn
&\!\!\!=\!\!\!& \tr[(\Pi_P-\Pi_N)\frac{d}{dt} X(t)]=2\ \tr[\Pi_P\frac{d}{dt} X(t)]\\
\nn
&\!\!\!=\!\!\!& 2\  \tr[-\Pi_P {\mathscr M}(X(t))]\\
\nn
&\!\!\!=\!\!\!& \!\! \int \!\! \tr   [\Pi_P[2LX(t)L^\ast -X(t)L^\ast L-L^\ast LX(t)]].\\
\label{33}
\eea
Define $X_P(t):=\Pi_P X(t)$,  and $X_N(t):=\Pi_N X(t)$. Notice that 
\beq\label{34}
X_P(t)=X_P(t)\Pi_P= \Pi_P X_P(t),
\eeq
and a similar equality holds for $X_N(t)$.
We also have $X_P(t)-X_N(t)=X(t)$ so that, by linearity, 
 the integrand   $\cal I$ of the last member of (\ref{33}) may be written as 
${\cal I}={\cal I}_P+{\cal I}_N$ with
\bea\nn
{\cal I}_P&\!\!\!:=\!\!\!&\tr\left[\Pi_P[2LX_P(t)L^\ast -X_P(t)L^\ast L-L^\ast LX_P(t)]\right]\\
\nn
&\!\!\!=\!\!\!&\tr[2\Pi_PLX_P(t)L^\ast -X_P(t)L^\ast L-\Pi_PL^\ast LX_P(t))]\\
\nn
&\!\!\!=\!\!\!&2\ \tr[L^\ast\Pi_P LX_P(t) - L^\ast L X_P(t))]\\
\nn
&\!\!\!=\!\!\!& 2\ \tr[L^\ast(\Pi_P - I) L X_P(t)]\\
&\!\!\!=\!\!\!& 2\ \tr[[X_P(t)]^{1/2}L^\ast(\Pi_P - I) L [X_P(t)]^{1/2}]\leq 0,
\eea
where we have used the cyclic  property of the trace operator, 
equality (\ref{34}), the fact that $X_P(t)={\rm e}^{-\alpha t}Y_P\geq 0$
and eventually that $\Pi_P - I\leq 0$ because $\Pi_P$ is an orthogonal projection.
As for ${\cal I}_N$, we have
\bea\nn
{\cal I}_N&\!\!\!:=\!\!\!&\tr\left[\Pi_P[-2LX_N(t)L^\ast +X_N(t)L^\ast L+L^\ast LX_N(t)]\right]\\
&\!\!\!=\!\!\!&-2\ \tr[\Pi_PLX_N(t)L^\ast\Pi_P]\leq 0,
\eea
where we have used the fact that $\Pi_P X_N(t)=0$, the cyclic property of the trace operator, the fact that $\Pi_P=\Pi_P^2$ and eventually that $X_N(t)\geq 0$.
In conclusion we have that ${\cal I}\leq 0$ so that
$\frac{d}{dt} \|X(t)\|_1=\int{\cal I}\leq 0$.
On the other hand 
\beq
0\geq\frac{d}{dt} \|X(t)\|_1=\frac{d}{dt} \|{\rm e}^{-\alpha t}Y\|_1=
-\alpha{\rm e}^{-\alpha t}\|Y\|_1.
\eeq
so that $\alpha\geq 0$ (recall that,  $Y$ being an eigenmatrix, it is not the zero matrix and hence  $\|Y\|_1>0$).
\cvd

Theorem \ref{theo33} and (\ref{sigb}) allow us to conclude that the matrix $B$ in (\ref{matrep}) is such that $\sigma(B)\subset [0, 1)$ and hence it is a discrete-time stability matrix.

\subsection{Center Manifold theory}

Let us go back to the original non-linear map $\Theta$.
The iteration (\ref{algori}) may be incrementally represented as
$\Lambda_{k+1}-\Lambdaopt=\Theta(\Lambda_{k})-\Lambdaopt$ and, by 
Taylor series expansion, as
\beq\label{incrementalev}
X_{k+1}={\mathscr M}(X_k) + m(X_k)
\eeq
where we have defined
$X_{k}:=\Lambda_{k}-\Lambdaopt$ and $m$ is the residue function that vanishes with its first derivatives at the origin.
Moreover, notice that, from (\ref{serveperprop}), (\ref{miden}), and (\ref{incrementalev}), we immediately get 
\beq\label{m=0}
 m(X_\perp)=0\quad \forall\ X_\perp\in(\Rgamma)^\perp.
\eeq
\begin{theo}
The set ${\cal L}_{\circ+}$ is locally asymptotically stable for $\Theta$.
\end{theo}
\IEEEproof
We resort again  to the coordinate representation of $X$ introduced in Remark \ref{remvectoriz}.
Moreover, we partition $x$ in the form  $x=\bmat{c}x^\perp\\x^r\emat$, where $x^\perp\in\C^{n_\perp}$ is the component of $x$ corresponding to $(\Rgamma)^\perp$, and
$x^r\in\C^{n_r}$. In these coordinates  the incremental evolution (\ref{incrementalev})
is represented by  
\beq\label{eqincrem}
\left\{
\begin{array}{l}
x^\perp_{k+1}=x^\perp_{k}+ f(x^\perp_{k},x^r_{k})\\
x^r_{k+1}=B x^r_{k}+ g(x^\perp_{k},x^r_{k}),
\end{array}\right.
\eeq
where, as already discussed, $B$ is a stability matrix.
We are now in the setting of Center Manifold theory, see \cite[pages 34--35]{CARR}.
The first and, in general, most difficult step to apply this theory is to find a center manifold, i.e. a ${\cal C}^2$ function $h:\C^{n_\perp}\rightarrow\C^{n_r}$ that vanishes with its first derivatives at the origin, and such that the center manifold equation
\beq\label{cmeq}
h(x^\perp+ f(x^\perp,h(x^\perp))=Bh(x^\perp)+g(x^\perp,h(x^\perp))
\eeq
is satisfied.
In our situation, however, this equation admits a solution that may be computed very
easily. In fact, in view of (\ref{m=0}), it is immediate to check that
\beq
f(x^\perp,0)=0,\ g(x^\perp,0)=0,\quad \forall x^\perp.
\eeq
Therefore, we may choose as a solution to equation (\ref{cmeq}) the  identically zero 
function $h$.
The asymptotic behavior of trajectories of (\ref{eqincrem}) originating in a neighborhood of the origin is determined by the flow on the center manifold whose dynamics is governed by the equation
\beq\label{cmfl}
u_{k+1}=u_k +f(u_k,h(u_k))=u_k +f(u_k,0)=u_k.
\eeq
Clearly, the zero solution of (\ref{cmfl}) is stable and thus, as stated in \cite[Theorem 8, page 35]{CARR}:
\begin{enumerate}
\item The zero solution of (\ref{eqincrem}) is stable.
\item
There exists a solution $u_k=u$ of (\ref{cmfl}) and two positive constants $\kappa$ and $\beta<1$, such that 
\beq
|x^\perp_{k}-u|\leq \kappa\beta^k,\quad |x^r_{k}|=|x^r_{k}-h(u)|\leq \kappa\beta^k.
\eeq
\end{enumerate}
In conclusion, if the initial condition of (\ref{eqincrem}) is sufficiently close to the origin, than 
$$x_k=\bmat{c}x_k^\perp\\x_k^r\emat\longrightarrow\bmat{c}\bar{x}^\perp\\0\emat,$$
i.e. $x_k$ converges to a state representing an element of $(\Rgamma)^\perp$.
This is equivalent to say that for any $\Lambdaopt\in{\cal L}_{\circ+}$ (defined in (\ref{callzp})) there exists a neighborhood ${\cal B}(\Lambdaopt)$ such that all trajectories $\{\Lambda_k\}$ generated by (\ref{algori}) and originating from ${\cal B}(\Lambdaopt)$, converge to ${\cal L}_{\circ+}$. Equivalently, 
${\cal L}_{\circ+}$ is locally asymptotically stable for $\Theta$.
\cvd

\section{Numerical implementation}\label{sec:NI}
In this section we discuss a numerically efficient implementation of the integral
in the iteration (\ref{algori}).
We show that it may be computed using very robust and reliable linear algebra algorithms.
We want to compute $\int G\frac{\Psi}{G^\ast\Lambda G}G^\ast$, 
where $G,\Psi$ and $\Lambda$ are given and $G^\ast\Lambda G$ is positive on $\T$.
To this aim, we assume that $\Psi$ is rational and, as a preliminary step we compute, using standard tools, a minimal minimum-phase spectral factor 
$$
W_\Psi(z) =H(zI-F)^{-1}G+D
$$
of $\Psi$.
We also employ the factorization
$$G^*\Lambda G = W^\ast(z) W(z),$$
with
\beq
W:=(B^{\ast}PB)^{-1/2}B^{\ast}PA(zI-A)^{-1}B+(B^{\ast}PB)^{1/2}
\eeq
derived in Appendix.
Thus, we clearly have
\beq\label{int56}
\int G\frac{\Psi}{G^\ast\Lambda G}G^\ast=\int (GW^{-1}W_\Psi)(GW^{-1}W_\Psi)^\ast
\eeq
so that the integral in (\ref{int56}) is the steady-state output covariance of the filter $GW^{-1}W_\Psi$ driven by normalized withe noise.
Then, let us compute a state space realization of $GW^{-1}W_\Psi$.
First, we observe that:
\beq
W^{-1}=[I-(B^{\ast}PB)^{-1}B^{\ast}PA(zI-Z)^{-1}B](B^{\ast}PB)^{-1/2},
\eeq
where $Z$, defined in (\ref{defz}), is a stability matrix.
Hence,
\bea\label{gwin}\nn
GW^{-1}&\!\!\!\!=\!\!\!\! &-(zI-A)^{-1}B(B^{\ast}PB)^{-1}B^{\ast}PA(zI-Z)^{-1}\times\\
\nn
&&B(B^{\ast}PB)^{-1/2}+(zI-A)^{-1}B(B^{\ast}PB)^{-1/2},\\
\eea
Notice that 
$$B(B^{\ast}PB)^{-1}B^{\ast}PA=A-Z=(zI-Z)-(zI-A).$$
Plugging this expression  into  (\ref{gwin}) we get 
\bea\nn
GW^{-1}&\!\!\!\!=\!\!\!\! &-(zI-A)^{-1}B(B^{\ast}PB)^{-1/2}\\
\nn
&&+(zI-A)^{-1}B(B^{\ast}PB)^{-1/2}\\
\nn
&&+(zI-Z)^{-1}B(B^{\ast}PB)^{-1/2}\\
&\!\!\!\!=\!\!\!\! &(zI-Z)^{-1}B(B^{\ast}PB)^{-1/2}.
\eea
Eventually, it is now easy to see that $GW^{-1}W_\Psi$ has the following state space realization
\beq
GW^{-1}W_\Psi=[0\mid I]\left(zI-\hat{F}\right)^{-1}\hat{G}
\eeq
with
$$
\hat{F}:=\bmat{cc}F&0\!\!\\ \!\!B(B^{\ast}PB)^{-1/2}H&Z\!\!\emat,\ \hat{G}:=\bmat{c}G\\\!\!B(B^{\ast}PB)^{-1/2}D\!\!\emat.$$
Notice that $\hat{F}$ is a stability matrix.
The following result comes now as a straightforward conclusion
\begin{prop}
Let $\Xi$ be the solution of the following discrete-time Lyapunov equation
\beq
\Xi=\hat{F}\Xi\hat{F}^\ast +\hat{G}\hat{G}^\ast.
\eeq
Then the integral in (\ref{int56}) is the bottom-right block of $\Xi$, i.e.,
\beq
\int G\frac{\Psi}{G^\ast\Lambda G}G^\ast=[0\mid I]\ \Xi\bmat{c}0
\\ I\emat.
\eeq
\end{prop}

In conclusion, for each iteration of the algorithm (\ref{algori}) we only have to compute: the solution of an algebraic Riccati equation of order $n$, the solution
of a discrete-time Lyapunov equation of order  $n+n_\Psi$ ($n_\Psi$ being the state space dimension of $W_\Psi(z)$), and the square root of a positive definite matrix $\Lambda$. All these operations are accomplished by standard linear algebra algorithms that may be implemented by numerically efficient and  robust routines.

\section{Evidence from simulations and a convergence conjecture}\label{sec:EFS}

In the previous section we have proven a local result.
In an extensive campaign  of simulations, however, we have always observed that the sequence $\{\Lambda_k\}$ converges very fast to a $\Lambdaopt$ in the closure $\overline{\cal L}_{\circ+}$ of ${\cal L}_{\circ+}$. Thus, we conjecture that indeed $\overline{\cal L}_{\circ+}$ is globally asymptotically stable for $\Theta$.
To this extent, we can do a few considerations. The function $\Theta$ maps the open ${\mathscr P}_+$ of positive definite matrices with  unitary trace to itself. Even if all fixed points in this open set are clearly in ${\cal L}_{\circ+}$, we cannot exclude that the sequence 
$\{\Lambda_k\}$ converges to the boundary of ${\mathscr P}_+$, i.e. to a singular matrix. Indeed, it is easy to see that there is a whole family of singular matrices in the boundary of ${\mathscr P}_+$  that are fixed points of $\Theta$. This is the family of the 
$1$-dimensional orthogonal projections. We have conducted some numerical experiments 
to understand the behavior of the map $\Theta$ in the neighborhood of $1$-dimensional orthogonal projections. We have observed that, even if we generate the sequence $\{\Lambda_k\}$ by choosing the initial condition arbitrarily close to a $1$-dimensional orthogonal projection $\Pi_1$, the sequence always converged to ${\cal L}_{\circ+}$.
For this reason we believe that, except for those in  $\overline{\cal L}_{\circ+}$, the $1$-dimensional orthogonal projections are {\em unstable} equilibrium points. A formal proof of this fact should   probably adopt a (nonlinear) Lyapunov approach. In fact, the derivative of the square root in the neighborhood of a singular matrix (as is an orthogonal projection) is infinite and thus a proof based on linearization does not seem  viable.

A second remark concerns the values of $\J_\Psi$ along the trajectory  $\{\Lambda_k\}$
generated by iterating $\Theta$. It has been shown in \cite{FERRANTE_PAVON_RAMPONI_FURTHERRESULTS} that
$\Delta \Lambda_k:= \Lambda_{k+1}-\Lambda_k$ is a descent direction for $\J_\Psi$.
Indeed, experimental evidence in numerical simulation is that more is true:
$\J_\Psi$ always decreases along trajectories $\{\Lambda_k\}$.
This fact, if proven, would be an important step toward a Lyapunov argument for {\em global convergence}.

\section{Conclusions}\label{sec:C}

The bottleneck of the spectral approximation techniques based on convex optimization is the solution of the dual problem \eqref{AP}.

So far, a closed form solution is not available for the general case, at least as long as the relative entropy functional is chosen as a spectral pseudo-distance. Thus,   an iterative algorithm has been proposed in  \cite{PAVON_FERRANTE_ONGEORGIOULINDQUIST,FERRANTE_PAVON_RAMPONI_FURTHERRESULTS} in  order to obtain the dual solution numerically. Necessary conditions for such an algorithm to be of interest are clearly its numerical efficiency and, most important, its convergence features.

While the proposed nonlinear iteration exhibits exceedingly good convergence features in all the simulation tests, a convergence result was missing. Our main result proves that the iteration is locally convergent to the manifold of full-rank solutions for the dual problem. The path to this results is quite tortuous yet provides new insights on the dynamics associated to the iteration, offering a potential standpoint towards the extension of the result to global convergence. 

We first develop a detailed analysis of the linearized map and of its spectral properties, also by resorting to an associated continuous-time evolution which provides us with a lower bound for the spectrum. Then the center manifold theorem, along with a desirable property of the nonlinear map, let us conclude that the manifold of solution is locally asymptotically stable.

The iterative algorithm presented in Section \ref{iteralg} is thus proven to be an eligible candidate   for being  the missing piece towards a  satisfactory, feasible solution of the spectral approximation problem in the general case.

As we conjectured in the previous section, supported by numerical simulations, we believe it should be possible to prove that convergence is almost global, namely that all the stationary points that are not in ${\mathcal L}_{0+}$ are in fact repulsive. Technical difficulties rule out a linearization approach, suggesting a general Lyapunov analysis as the natural pathway to the desired result. This indeed represents the most challenging yet compelling direction for further work.

We wish to thank prof.\ Christopher Byrnes for the insightful
discussions and for the advice he gave us during his staying in Padova.

\appendix
In the following we present a factorization result that is repeatedly used in different parts of the paper. 
\begin{lem}\label{lemapp}
Let
$G(z)=(zI-A)^{-1}B$ with $A\in \C^{n\times n},$ $B\in\C^{n\times 1},$ and $(A,B)$  a  reachable pair.
Let $\Lambda\in\Hermitian_n$ be such that
$G^*\Lambda G>0$ on $\T$.
Then, the following factorization holds:
\beq
G^*\Lambda G=W^*W,
\eeq
where
\beq\label{wlem}
W:=(B^{\ast}PB)^{-1/2}B^{\ast}PA(zI-A)^{-1}B+(B^{\ast}PB)^{1/2}
\eeq
and $P\in\Hermitian_n$ is the stabilizing solution of 
the algebraic Riccati equation
\beq\label{AREL}
\Pi =A^{\ast}\Pi A-A^{\ast}\Pi B(B^{\ast}\Pi B)^{-1}B^{\ast}\Pi A+\Lambda,
\eeq
so that the spectrum of closed loop matrix 
\beq\label{defz}
Z:=A-B(B^{\ast}PB)^{-1}B^{\ast}PA
\eeq
  lays inside the open unit disk.
\end{lem}
\IEEEproof
Given $G(z)=(zI-A)^{-1}B$ and an arbitrary
$\Pi\in\Hermitian_n$, the following identity holds \cite{FERRANTE_COLANERI_ALGEBRAICRICCATI}:
\beq\label{idenzero}
[G^*(z)\mid 1]\bmat{cc}A^{\ast}\Pi A-\Pi&A^{\ast}\Pi B\\B^{\ast}\Pi A&B^{\ast}\Pi B\emat \bmat{c}G(z)\\1\emat
\equiv 0.
\eeq
Therefore, we have
\beq\label{rgslg}
G^*\Lambda G=[G^*\mid 1]\bmat{cc}\Lambda+A^{\ast}\Pi A-\Pi &A^{\ast}\Pi B\\B^{\ast}\Pi A&B^{\ast}\Pi B\emat \bmat{c}G\\1\emat
\eeq
Since $G^\ast\Lambda G$ is positive on the whole $\T$,
there exists the stabilizing solution $P$ 
of the algebraic Riccati equation (\ref{AREL}).

Thus if we set $\Pi=P$,   the block matrix on  the right hand side of (\ref{rgslg}) has the  the following factorization
$$
\bmat{c}A^{\ast}PB\\
B^{\ast}PB\emat (B^{\ast}PB)^{-1} [B^{\ast}PA\mid B^{\ast}PB]$$
so that
\bea\nn
G^*\Lambda G&\!\!\!\!=\!\!\!\!&[G^*A^{\ast}PB+B^{\ast}PB] (B^{\ast}PB)^{-1} [B^{\ast}PAG+B^{\ast}PB]\\\label{last}
&\!\!\!\!=\!\!\!\! & W^\ast(z) W(z),
\eea
with $W(z)$ given by (\ref{wlem}).
\cvd

\bibliographystyle{plain}
\bibliography{fpr}

 \end{document}